\documentclass[12pt]{amsart}
\usepackage{amsfonts,amssymb,amsmath,amsthm}
\usepackage{url,mathabx}
\usepackage{enumerate}

\title[The general linear group as a complete invariant]{The general linear group as a complete invariant for C*-algebras}

\author{Thierry Giordano${}^1$}

\address{Department of Mathematics and Statistics\\
University of Ottawa\\
585 King Edward Avenue\\
Ottawa, Ontario\\
KiN 6N5\\
Canada}

\email{giordano@uottawa.ca}
\thanks{${}^1$T.G.\ was partially supported by a grant from NSERC Canada.}

\author{Adam Sierakowski${}^2$}

\address{School of Mathematics and Applied Statistics \\
Building 15\\
University of Wollongong\\
Wollongong NSW \\
2522\\
Australia}

\email{asierako@uow.edu.au}
\thanks{${}^2$Communicating author. A.S.\ was supported by the Australian Research Council.}

\theoremstyle{proclaim}
\newtheorem{theorem}{Theorem}[section]
\newtheorem{lemma}[theorem]{Lemma}
\newtheorem{corollary}[theorem]{Corollary}
\newtheorem{proposition}[theorem]{Proposition}
\theoremstyle{definition}
\newtheorem{remark}[theorem]{Remark}
\newtheorem{definition}[theorem]{Definition}

\newtheorem{notation}[theorem]{Notation}

\newcommand{\F}{\mathcal{F}}

\newcommand{\I}{\mathcal{I}}
\newcommand{\N}{\mathcal{N}}
\newcommand{\U}{\mathcal{U}}
\newcommand{\Z}{\mathcal{Z}}
\newcommand{\cP}{\mathcal{P}}
\newcommand{\cO}{\mathcal{O}}
\newcommand{\CC}{\mathbb{C}}
\newcommand{\NN}{\mathbb{N}}
\newcommand{\ZZ}{\mathbb{Z}}

\newcommand{\is}{\sim_c}
\newcommand{\esim}{\sim_s}
\newcommand{\hsim}{\sim_h}
\newcommand{\IA}{\widetilde{\I(A)}}
\newcommand{\IB}{\I(B)\widetilde{\ }}
\newcommand{\PA}{\widetilde{\cP(A)}}
\newcommand{\PB}{\widetilde{\cP(B)}}

\newcommand{\ie}{\textit{i.e.},~}
\newcommand{\cst}{$\rm{C}^*$}

\def \bib(#1;#2;#3;#4;#5;#6) {{#1}, {\it #2} {#3},
{\bf#4} (#5) {#6}\par\smallskip}

\date{\today}

\numberwithin{equation}{section}

\begin{document}

\begin{abstract}
In 1955 Dye proved that two von Neumann factors not of type $I_{2n}$ are isomorphic if and only if their unitary groups are isomorphic as abstract groups. We consider an analogue for \cst-algebras and show that the topological general linear group is a classifying invariant for simple unital AH-algebras of slow dimension growth and of real rank zero, and that the abstract general linear group is a classifying invariant for unital Kirchberg algebras in the UCT class.
\end{abstract}


\maketitle

\section*{Introduction}
Since the introduction of the Elliott invariant as a classifying invariant for \cst-algebras, the classification program for \cst-algebras has been rapidly evolving. New invariants were introduced to enrich the program, some more general, and other tailored to specific applications. For a large class of simple, amenable, unital, separable \cst-algebras, Al-Rawashdeh, Booth and the first named author showed in \cite{AlRBooGio} that their unitary group forms a classifying invariant: From an isomorphism of the unitary group of such algebras, they deduced an isomorphism of their Elliott invariant. 

In this paper we look at the general linear group (\ie{}the group of invertible elements) of unital \cst-algebras as an invariant.  For each unital \cst-algebra $A$ we will denote its general linear group by $GL(A)$ and its set of idempotents by $\I(A)$ (see Notation \ref{note}). Given two unital \cst-algebras $A$ and $B$, and a group isomorphism $\varphi\colon GL(A)\to GL(B)$ between their general linear groups, the formula
\begin{equation*}
1-2\theta_\varphi(e)=\varphi(1-2e), \ \ \ e\in \I(A),
\end{equation*}
induces a bijection $\theta_\varphi\colon \I(A)\to \I(B)$ between the set of idempotents of $A$ and $B$. This map is not in general an orthoisomorphism of idempotents (\ie{}a bijective map which preserves
orthogonality of commuting idempotents). However, it turns out that in many cases $\theta_\varphi$ is essentially an orthoisomorphism. More precisely, generalising the notion of oddly decomposability given in \cite{AlRBooGio} (see Definition~\ref{oddly}), we show in Theorem~\ref{thm2} that there exist a partitioning of the non-trivial elements of $\I(A)$ into two set $\I_o,\I_{\bar{o}}$, such that the map  $\tilde\theta_\varphi  \colon \I(A)\to \I(B)$ defined by
\begin{align*}
\tilde\theta_\varphi (e)=\left\{
\begin{array}{ll}
\theta_\varphi (e),&\mbox{ if } e \in \I_o\\
1-\theta_\varphi (e),&\mbox{ if } e \in \I_{\bar{o}}\\
1,&\mbox{ if } e =1\\
0,&\mbox{ if } e =0
\end{array}\right.
\end{align*}
is an orthoisomorphism. Using the maps $\tilde\theta_\varphi$ and $\varphi$ between the idempotents and invertibles of $A$ and  $B$ we construct homomorphisms from $K_0(A)$ to $K_0(B)$ and from $K_1(A)$ to $K_1(B)$ and invoke on classification to show $A$  and $B$ are isomorphic. By investigating which \cst-algebras are oddly decomposable we prove the following two main results
\begin{itemize}
\item[(i)] Let $A$ and $B$ be simple, unital $AH$-algebras of slow dimension growth and of real rank zero. Then $A$ and $B$ are isomorphic if and only if their general linear groups are topologically isomorphic.
\item[(ii)] Let $A$ and $B$ be unital Kirchberg algebras in the UCT class. Then $A$ and $B$ are isomorphic if and only if their general linear groups are isomorphic as abstract groups.
\end{itemize}
In the case the algebras $A$ and $B$ are simple and finite dimensional we refer to \cite{SchWar} by Schreier and Van der Waerden (see also \cite{GioSie3}, \cite{HouHua}, and \cite{Sem} for related results).

\section{Properties of the induced map $\theta_\varphi$}\label{sec2}
Let $A$ and $B$ be two unital \cst-algebras. If $\varphi\colon GL(A) \to GL(B)$ is a group homomorphism between the general linear groups of $A$ and $B$, then $\varphi$ defines a map $\theta=\theta_\varphi  \colon \I(A)\to \I(B)$ by setting
\begin{equation*}
1-2\theta_\varphi(e)=\varphi(1-2e), \ \ \ e\in \I(A).
\end{equation*}
A simple computation shows that $\theta_\varphi(e)$ is an idempotent for each $e\in \I(A)$ making $\theta_\varphi$ well defined. If $\varphi$ is moreover a bijection---or more generally if $\varphi$ restricts to a bijection of symmetries (elements whose square equals the unit)---it follows that $\theta_\varphi$ is a bijection of idempotents. The following additional properties of the map $\theta$ can be easily checked by adapting the arguments in \cite{AlRBooGio} and \cite{Dye} to the present situation.

\begin{proposition}\label{prop1}
Let $A$ and $B$ be unital \cst-algebras, $\varphi\colon GL(A) \to GL(B)$ be a group isomorphism and $\theta$ be the induced map between idempotents. Then
\begin{enumerate}[(i)]
\item $\theta(ueu^{-1})=\varphi(u)\theta(e)\varphi(u)^{-1}$;
\item $\theta(0)=0$;
\item if $e,f\in \I(A)$ commute, then so do $\theta(e)$ and $\theta(f)$ in $\I(B)$;
\item $\theta(e\triangle f)=\theta(e)\triangle\theta(f)$, where $\triangle$ denotes the symmetric difference of commuting idempotents, \ie $e\triangle f = e+f-2ef$.
\end{enumerate}
\end{proposition}

If the center $\Z(B)$ of a unital \cst-algebra $B$ is reduced to the scalars, and $\varphi\colon GL(A) \to GL(B)$ is as above, then $\varphi(-1)=-1$: Indeed, note that $-1$ is a central element which is not $1$, but its product with itself equals $1$. The same is true for $\varphi(-1)$. As a consequence we get the following lemma, see also \cite{AlRBooGio} and \cite{Dye}.

\begin{lemma}\label{lem3}
Let $A$ and $B$ be unital \cst-algebras, whose center $\Z(B)=\CC 1$. Let $\varphi\colon GL(A) \to GL(B)$ be a group isomorphism and $\theta\colon \I(A)\to \I(B)$ be as above. Then $\theta(1) = 1$, and for each $e \in \I(A)$, $\theta(1-e) = 1-\theta(e)$.
\end{lemma}

To simplify notation, let us introduce the following:

\begin{notation}\label{note}
(i) The quadruple $(A,B,\varphi,\theta)$ will denote a pair of simple unital \cst-algebras $A$ and $B$, a group isomorphism $\varphi\colon GL(A) \to GL(B)$, and the induced bijection $\theta \colon \I(A)\to \I(B)$. 

(ii) Let $A$ be a unital \cst-algebra. Denote by $\I(A)$ the set of idempotents in $A$, and by $\IA$ the set $\I(A)\backslash\{0,1\}$ of \emph{non-trivial} idempotents in $A$. Denote by $GL(A)$ the general linear group of invertible elements in $A$.
\end{notation}

\begin{definition}
Let $A$ be a unital \cst-algebra. We say that two idempotents $e,f\in \I(A)$ are \emph{similar}, denoted $e\esim f$, if there exist $u\in GL(A)$ such that
\begin{equation*}
f=ueu^{-1}.
\end{equation*}
\end{definition}

The following lemma is a generalisation of Lemma~10 in \cite{Dye} to simple, unital \cst-algebras.

\begin{lemma}\label{lem1}
Let $(A,B,\varphi,\theta)$ be as in (\ref{note}). Then for each fixed $e\in \I(A)$,
\begin{equation*}
\varphi(\lambda e +1 - e) \in \CC\theta(e)+\CC\theta(1-e), \ \ \ \lambda\in \CC\backslash\{0\}.
\end{equation*}
\end{lemma}

\begin{proof}
Fix $\lambda\in \CC\backslash\{0\}$ and set $x:= \varphi(\lambda e +1 - e)$. Since every idempotent in $B$ is similar to projection (by Proposition~4.6.2 in \cite{Bla}) we can choose $u\in GL(B)$ such that $q=u\theta(e)u^{-1}$ is a projection. For any subset $S \subseteq B$, let $S^{\prime}$ denote its relative commutant in $B$, and $S^{''}$ its (relative) bicommutant.

We show that $uxu^{-1}\in \{q\}'$. Since $q=u\theta(e)u^{-1}$ we just need to show that $x\theta(e)=\theta(e)x$. This follows from $x\varphi(1-2e)=\varphi(1-2e)x$.

We show $uxu^{-1}\in \{q\}''$. Fix $b\in \{q\}'$. Since $q$ is selfadjoint $\{q\}'$ is a \cst-subalgebra of $B$ and contains unitary elements $\varphi(u_1),\dots,\varphi(u_4)$ that span $b$, for some $u_1,\dots, u_4\in GL(A)$. Using $\varphi(u_i)\in \{q\}'$ commutes with $q=u\theta(e)u^{-1}$ we have that $u^{-1}\varphi(u_i)u$ commutes with $\theta(e)$ and with $\varphi(1-2e)$. This implies that $\varphi^{-1}(u^{-1})u_i\varphi^{-1}(u)$ commutes with $1-2e$, with $e$ and with $\lambda e +1 - e$. We now have that  $u^{-1}\varphi(u_i)u$ commutes with $x=\varphi(\lambda e +1 - e)$. Therefore $\varphi(u_i)$ commutes with $uxu^{-1}$, and $b$ commutes with $uxu^{-1}$. Since $b$ was an arbitrary element in $\{q\}'$ we conclude that $uxu^{-1}$ commutes with every element in $\{q\}'$, \ie{}$uxu^{-1}\in \{q\}''$.

Since $q$ is a projection $\{q\}''\cap \{q\}'=\CC q + \CC (1-q)$, using the fact that $B$ is simple so the hereditary \cst-subalgebra $qBq$ is simple and consequently has centre $\CC q$ (similarly for $1-q$). Multiplying $uxu^{-1}$ on the left by $u^{-1}$ and on the right by $u$ we see that $x\in \CC\theta(e)+\CC\theta(1-e)$.
\end{proof}

In the following we let $\CC^*$ denote the group $(\CC\backslash\{0\}, \cdot)$ of non-zero complex numbers with multiplication as the group operation.

\begin{lemma}\label{lem2}
Let $(A,B,\varphi,\theta)$ be as in (\ref{note}). Then for each fixed $e\in \IA$ there exist group homomorphisms $a_e,b_e\colon \CC^*\to \CC^*$ such that
\begin{equation*}
\varphi(\lambda e +1 - e)= a_e(\lambda)\theta(e)+b_e(\lambda)\theta(1-e), \ \ \ \lambda\in \CC\backslash\{0\}.
\end{equation*}
\end{lemma}
\begin{proof}
Fix $\lambda\in \CC\backslash\{0\}$. Since $e\in \IA$ the elements $\theta(e)$ and $\theta(1-e)$ are nonzero, and by Lemma~\ref{lem3} they are linearly independent. Using Lemma~\ref{lem1} we therefore have unique coefficients $a,b\in\CC$ such that
\begin{equation*}
\varphi(\lambda e +1 - e) = a\theta(e)+b\theta(1-e).
\end{equation*}
Assuming $b=0$ we obtain $\varphi(\lambda e +1 - e)^2=a^2\theta(e)=a(a\theta(e))=a\varphi(\lambda e +1 - e)$. Hence $a\theta(e)=\varphi(\lambda e +1 - e)=a1=a\theta(1)$. Since $\varphi(\lambda e +1 - e)$ is invertible we get $\theta(e)=\theta(1)$, contradiction contradicting injectivity of $\theta$. By symmetry both $a,b\in \CC\backslash\{0\}$. 

It is easy to see that $a_e, b_e$ are multiplicative and unital using that the map $\lambda\mapsto \varphi(\lambda e +1 - e)$ is multiplicative and unital. We conclude both maps are group homomorphisms.
\end{proof}

Since the maps $a_e, b_e$ are group homomorphism of $\CC^*$ we will use their (multiplicative) inverses without any further explanation. To each $e\in \IA$, we associate the pair of maps $(a_e,b_e)$ and the group homomorphism $c_e:=a_eb_e^{-1}$ of $\CC^*$. Moreover, we denote by $\is$ the equivalence relation on $\IA$, given by: 
\begin{equation*}
e\is f \ \ \ \textrm { if and only if } \ \ \ c_e=c_f.
\end{equation*}

The following proposition is essentially Proposition 2.7 in \cite{AlRBooGio}. The proof can be adapted to the present situation and it is left to the reader.
\begin{proposition}\label{prop2}
Let $(A,B,\varphi,\theta)$ be as in (\ref{note}). Then for each $e\in \IA$
\begin{enumerate}[(i)]
\item If $f \in \IA$ is similar to $e$ then $e \is f$;
\item $c_e(\lambda)^2\neq 1$, for every $\lambda\in \CC\backslash\{-1,0,1\}$;
\item $c_e=c_{1-e}$.
\end{enumerate}
\end{proposition}

\begin{definition}
Two or more idempotents in a \cst-algebra $A$ are {\em orthogonal} provided that any two of these idempotents commute and their product is equal to zero.
\end{definition}

\begin{lemma}\label{lem4}
Let $(A,B,\varphi,\theta)$ be as in (\ref{note}). Suppose that $e,f\in \IA$ are two orthogonal idempotents in $A$. Then
\begin{align*}
\theta(e+f)&=\theta(e)\theta(1-f)+\theta(1-e)\theta(f)\\
\theta(1-e-f)&=\theta(e)\theta(f)+\theta(1-e)\theta(1-f)
\end{align*}
\end{lemma}

\begin{proof}
The proof of Lemma~2.3 in \cite{AlRBooGio} does not directly generalise, so we include a short proof: Using Proposition~\ref{prop1} and Lemma~\ref{lem3} we have
\begin{equation*}
\theta(e+f)=\theta(e)\triangle\theta(f)=\theta(e)+\theta(f)-2\theta(e)\theta(f)=\theta(e)\theta(1-f)+\theta(1-e)\theta(f).
\end{equation*}
The second equality follows by subtracting both sides of the above equality from 1.
\end{proof}

\begin{remark}
The proof of Theorem~\ref{thm1}, Corollary~\ref{cor1}, and Corollary~\ref{cor2} corresponds to Proposition~2.8 and Theorem~2.9 in \cite{AlRBooGio}, but our proof includes a new characterisation of when $c_e=c_f$, $\dots$, $c_e=c_{e+f}^{-1}$ in terms of the equations (\ref{eq5})-(\ref{eq8}). This observation is essential in the subsequent proofs of Corollary~\ref{cor1}, Corollary~\ref{cor2} used to prove Lemma~\ref{lem5} and Lemma~\ref{lem6}.
\end{remark}

\begin{theorem}\label{thm1}
Let $(A,B,\varphi,\theta)$ be as in (\ref{note}).  Suppose that $e,f\in \IA$ are two orthogonal idempotents in $A$ not adding to one. Then
\begin{align*}
\theta(e)\theta(f)=0 \ \ \ &\Leftrightarrow \ \ \ c_e=c_f=c_{e+f}\\
\theta(1-e)\theta(1-f)=0\ \ \ &\Leftrightarrow \ \ \  c_e=c_f=c_{e+f}^{-1}\\
\theta(1-e)\theta(f)=0\ \ \ &\Leftrightarrow \ \ \ c_e=c_f^{-1}=c_{e+f}\\
\theta(e)\theta(1-f)=0\ \ \ &\Leftrightarrow \ \ \ c_e=c_f^{-1}=c_{e+f}^{-1}
\end{align*}
\end{theorem}
\begin{proof}
Since $\varphi(\lambda e+1-e)\varphi(\lambda f+1-f)=\varphi(\lambda (e+f)+1-(e+f))$ for $\lambda\neq 0$, Lemma~\ref{lem2} ensures that \begin{equation*}
(a_e\theta(e)+b_e\theta(1-e))(a_f\theta(f)+b_f\theta(1-f))=a_{e+f}\theta(e+f)+b_{e+f}\theta(1-(e+f)).
\end{equation*}
Multiplying each side of the above equality by $\theta(e)\theta(f)$, $\theta(1-e)\theta(f)$, $\theta(1-e)\theta(f)$ or $\theta(e)\theta(1-f)$ and then using Lemma~\ref{lem4}, we obtain the following four equations
\begin{align}
a_ea_f\theta(e)\theta(f)&=b_{e+f}\theta(e)\theta(f)\label{eq1}\\
b_eb_f\theta(1-e)\theta(1-f)&=b_{e+f}\theta(1-e)\theta(1-f)\label{eq2}\\
b_ea_f\theta(1-e)\theta(f)&=a_{e+f}\theta(1-e)\theta(f)\label{eq3}\\
a_eb_f\theta(e)\theta(1-f)&=a_{e+f}\theta(e)\theta(1-f)\label{eq4}
\end{align}
Consider the following properties
\begin{align}
a_ea_f &=b_{e+f}\label{eq5}\\
b_eb_f &=b_{e+f}\label{eq6}\\
b_ea_f &=a_{e+f}\label{eq7}\\
a_eb_f &=a_{e+f}\label{eq8}
\end{align}
We claim that

\begin{align*}
(\ref{eq7}), (\ref{eq8})\ \ \ &\Leftrightarrow  \ \ \ c_e=c_f\\
(\ref{eq5}), (\ref{eq6})\ \ \ &\Leftrightarrow \ \ \  c_e=c_f^{-1}\\
(\ref{eq6}), (\ref{eq8})\ \ \ &\Leftrightarrow \ \ \ c_e=c_{e+f}\\
(\ref{eq5}), (\ref{eq7})\ \ \ &\Leftrightarrow \ \ \ c_e=c_{e+f}^{-1}
\end{align*}
Going from left to right is straight forward. To go from right to left one simply adds two of the equations (\ref{eq1})-(\ref{eq4}), possibly with coefficients. For example if $c_e=c_f$ then $a_eb_f=b_ea_f$. By Lemma~\ref{lem4} the equality $(\ref{eq3})+(\ref{eq4})$, where we add each side separately, reduces to 
\begin{equation*}
a_eb_f\theta(e+f)=b_ea_f\theta(e+f)=a_{e+f}\theta(e+f).
\end{equation*}
Hence (\ref{eq7}), and (\ref{eq8}) both hold. The remaining three equivalences are obtained similarly using $(\ref{eq1})+(\ref{eq2})$, $a_e\cdot (\ref{eq2})+ b_e\cdot (\ref{eq3})$, and $b_e\cdot (\ref{eq1})+ a_e\cdot (\ref{eq3})$. We obtain
\begin{align*}
(\ref{eq6}), (\ref{eq7}),(\ref{eq8}) \ \ \ &\Leftrightarrow \ \ \ c_e=c_f=c_{e+f}\\
(\ref{eq5}), (\ref{eq7}),(\ref{eq8})\ \ \ &\Leftrightarrow \ \ \  c_e=c_f=c_{e+f}^{-1}\\
(\ref{eq5}), (\ref{eq6}),(\ref{eq8})\ \ \ &\Leftrightarrow \ \ \ c_e=c_f^{-1}=c_{e+f}\\
(\ref{eq5}), (\ref{eq6}),(\ref{eq7})\ \ \ &\Leftrightarrow \ \ \ c_e=c_f^{-1}=c_{e+f}^{-1}
\end{align*}
We now show that $\theta(e)\theta(f)=0$ if and only if $c_e=c_f=c_{e+f}$. The other three equivalences follow from similar calculations. 

Suppose first that ${\theta(e)\theta(f)=0}$. Assume that $c_e=c_f=c_{e+f}$ does not hold. We derive a contradiction. Using the observation above, one of (\ref{eq6}), (\ref{eq7}), or (\ref{eq8}) does not hold. It follows that $\theta(1-e)\theta(1-e)=0$, $\theta(1-e)\theta(f)=0$, or $\theta(e)\theta(1-f)=0$. Adding this to $\theta(e)\theta(f)=0$, Lemma~\ref{lem4} gives $\theta(1-e-f)=0$, $\theta(f)=0$, or $\theta(e)=0$. Contradiction.

Conversely suppose that the first equation $c_e=c_f=c_{e+f}$ above holds. Since $c_g^2\neq 1$ for $g=e,f,e+f$, by Proposition~\ref{prop2}$(ii)$, we obtain that all the other three equations are false. Since (\ref{eq6}), (\ref{eq7}), and (\ref{eq8}) hold, (\ref{eq5}) must fail. We conclude that $\theta(e)\theta(f)=0$.
\end{proof}

\begin{corollary}\label{cor1}
Let $(A,B,\varphi,\theta)$ be as in (\ref{note}). Suppose that $e,f\in \IA$ are two orthogonal, $\is$-equivalent, idempotents in $A$ not adding to one. Then precisely one of $\theta(e)\theta(f)$, $\theta(1-e)\theta(1-f)$ is zero.
\end{corollary}
\begin{proof}
If both terms are zero then $c_{e+f}=c_{e+f}^{-1}$. If both terms are non-zero then (\ref{eq5}) and (\ref{eq6}) are true: If (\ref{eq5}) fails then $\theta(e)\theta(f)=0$ by (\ref{eq1}), and similarly for (\ref{eq6}). Hence $c_e=c_f^{-1}$ and by $\is$-equivalence also $c_{f}=c_{f}^{-1}$.
\end{proof}
 
\begin{corollary}\label{cor2}
Let $(A,B,\varphi,\theta)$ be as in (\ref{note}). Suppose that $e,f\in \IA$ are two orthogonal idempotents in $A$ not adding to one. Then precisely one of $c_{e}=c_{e+f}$, $c_{e}=c_{e+f}^{-1}$ is true.
\end{corollary}
\begin{proof}
Consider the four equations from Theorem~\ref{thm1}
\begin{align*}
&c_{e}=c_{f}=c_{e+f} &c_{e}=c_{f}=c_{e+f}^{-1}\\
&c_{e}=c_{f}^{-1}=c_{e+f} &c_{e}=c_{f}^{-1}=c_{e+f}^{-1}
\end{align*}
It is clear that at least one of (\ref{eq5})-(\ref{eq8}) is false.  Consequently at least one of the four equations above is true: If (\ref{eq5}) fails then $\theta(e)\theta(f)=0$, so $c_{e}=c_{f}=c_{e+f}$. Similarly for the remaining cases. We obtain $c_{e}=c_{e+f}$ or $c_{e}=c_{e+f}^{-1}$ is true. Both can not be true because as this contradicts Proposition~\ref{prop2}$(ii)$.
\end{proof}

\begin{definition}
Let $(A,B,\varphi,\theta)$ be as in (\ref{note}). For any subset $S\subseteq \I(A)$ we say that $\theta\colon\I(A)\to \I(B)$ \emph{preserves orthogonality} (resp.\ \emph{flips orthogonality}) on $S$ if $\theta(e)\theta(f)=0$ (resp.\ $\theta(1-e)\theta(1-f)=0$) for any two orthogonal idempotents $e$ and $f$ in $S$.
\end{definition}

\begin{lemma}\label{lem5}
Let $(A,B,\varphi,\theta)$ be as in (\ref{note}). Suppose that $e,f,g\in \IA$ are three orthogonal, $\is$-equivalent, idempotents in $A$ not adding to one. If $\theta$ preserves (resp.\ flips) orthogonality on a subset of $\{e,f,g\}$ of size two, then $\theta$ preserves (resp.\ flips) orthogonality on all of $\{e,f,g\}$.
\end{lemma}

\begin{proof}
The proof of Lemma~2.14 in \cite{AlRBooGio} does not generalise nicely, so we include a short proof:
	
(1) Suppose that $\theta(e)\theta(f)=0$, $\theta(e)\theta(g)=0$. Assume for contradiction that $\theta(1-f)\theta(1-g)=0$. By Lemma~\ref{lem4}
\begin{equation*}
\theta(e)\theta(f+g)=\theta(e)\big(\theta(f)\theta(1-g)+\theta(1-f)\theta(g)\big)=0.
\end{equation*} 
Theorem~\ref{thm1} implies that $c_e=c_{f+g}=c_{e+f+g}$ and $c_f=c_g=c_{f+g}^{-1}$. Hence $c_{f+g}=c_{f+g}^{-1}$ which is false. We get $\theta(1-f)\theta(1-g)\neq 0$. By Corollary~\ref{cor1} we conclude $\theta(f)\theta(g)=0$.

(2) Now suppose that $\theta(1-e)\theta(1-f)=0$, $\theta(1-e)\theta(1-g)=0$. Assume for contradiction that $\theta(f)\theta(g)=0$. By Lemma~\ref{lem4},
\begin{equation*}
\theta(1-e)\theta(f+g)=\theta(1-e)\big(\theta(f)\theta(1-g)+\theta(1-f)\theta(g)\big)=0.
\end{equation*} 
Theorem~\ref{thm1} implies that $c_e=c_{f+g}^{-1}=c_{e+f+g}$ and $c_f=c_g=c_{f+g}$. We obtain $c_{f+g}=c_{f+g}^{-1}$ which is false. Consequently, we have that $\theta(f)\theta(g)\neq 0$. By Corollary~\ref{cor1}, $\theta(1-f)\theta(1-g)=0$. It is evident (1)-(2) suffice to complete the proof.
\end{proof}

\begin{lemma}\label{lem6}
Let $(A,B,\varphi,\theta)$ be as in (\ref{note}). Suppose that $e,f,g\in \IA$ are three orthogonal, $\is$-equivalent, idempotents in $A$ not adding to one. Then 
\begin{equation*}
e\is f\is g\is e+f+g.
\end{equation*}
\end{lemma}
\begin{proof}
The proof of Lemma~2.16 in \cite{AlRBooGio} can be adapted to the present situation. We include the details for completeness: By Corollary~\ref{cor1} precisely one of $\theta(f)\theta(g)$, $\theta(1-f)\theta(1-g)$ is zero. If $\theta(f)\theta(g)=0$, then by Lemma~\ref{lem4} and Lemma~\ref{lem5} 
\begin{equation*}
\theta(e)\theta(f+g)=\theta(e)\big(\theta(f)\theta(1-g)+\theta(1-f)\theta(g)\big)=0.
\end{equation*}
Using Theorem~\ref{thm1}, $c_e=c_{f+g}=c_{e+f+g}$. Similarly, if $\theta(1-f)\theta(1-g)=0$, then by Lemma~\ref{lem4} and Lemma~\ref{lem5} 
\begin{equation*}
\theta(1-e)\theta(f+g)=\theta(1-e)\big(\theta(f)\theta(1-g)+\theta(1-f)\theta(g)\big)=0.
\end{equation*}
Using Theorem~\ref{thm1}, $c_e=c_{f+g}^{-1}=c_{e+f+g}$.
\end{proof}

\section{Oddly decomposable \cst-algebras}
\label{sec3}
Let $(A,B,\varphi,\theta)$ be as in (\ref{note}), and let $\is$ be the equivalence relation on $\IA$ introduced in Section \ref{sec2}. We now introduce a sufficient condition on the \cst-algebra $A$, such that $\IA/\negthickspace\is$ has at most two elements.
 
\begin{definition}
\label{oddly}
A unital \cst-algebra $A$ is said to be \emph{oddly decomposable} if for every pair of idempotents $e,f\in \IA$ there is an odd integer $n\geq 3$ and $n$ orthogonal idempotents $g_1,\dots,g_n\in \IA$ adding to $f$, such that each $g_i$ is similar to some $g_i'\in\IA$ with $g_i'e=eg_i'=g_i'\neq e$.
\end{definition}

\begin{remark}
\label{rem3.1}
Oddly decomposable \cst-algebras where introduced in \cite{AlRBooGio}, but with a definition in terms of projections and unitary equivalence rather what idempotents and similarity. In \cite{AlRBooGio} a unital \cst-algebra $A$ was called oddly decomposable if for every pair of projections $p,q \in A\backslash\{0,1\}$ there exist an odd integer $n\geq 3$ and $n$ orthogonal non-zero projections $r_1,\dots,r_n\in A$ adding to $q$, such that each $r_i$ is unitary equivalent to some projection $r_i'\in A$ with $r_i'< p$. Let us outline why the two definitions coincide:

Fix a pair of idempotents $e,f\in \IA$. Find projections $p,q\in A$ and invertible elements $u,v\in GL(A)$ such that $e=upu^{-1}$, and $f=vqv^{-1}$ (see Lemma~\ref{lem11}). Assuming odd decomposability in sense of \cite{AlRBooGio} there exist an odd integer $n\geq 3$ and a decomposition of $q$ as a sum $q = \sum_{i=1}^n r_i$ of pairwise nonzero orthogonal projections $r_i$ of $A$, such that each $r_i$ is unitarily equivalent to some projection $r_i' < p$. Define
\begin{equation*}
g_i:=vr_iv^{-1}, \ \ \ g_i':=ur_i'u^{-1}, \ \ \ i=1,\dots,n.
\end{equation*}
It follows that $g_i,g_i'\in \IA$ have the properties needed to make $A$ oddly decomposable in sense of Definition~\ref{oddly}.

Conversely, fix a pair of projections $p,q \in A\backslash\{0,1\}$. Assuming odd decomposability in sense of Definition~\ref{oddly} there exist an odd integer $n\geq 3$ and $n$ orthogonal idempotents $g_1,\dots,g_n\in \IA$ adding to $q$, such that each $g_i$ is similar to some $g_i'\in\IA$ with $g_i'p=pg_i'=g_i'\neq p$. We can select $w_1,\dots, w_n\in GL(A)$ such that 
\begin{equation*}
r_i:=w_ig_iw_i^{-1}, \ \ \ i=1,\dots,n
\end{equation*}
are orthogonal projections adding to $q$. We can select projections $r_1',\dots, r_n'$ in $A$ such that each $r_i'$ is similar to $g_i'$ with $r_i'<p$. It follows that each $r_i$ is similar and hence also unitarily equivalent (see Lemma~\ref{lem11}) to $r_i'$, making $A$ oddly decomposable in sense of \cite{AlRBooGio}.
\end{remark}

\begin{notation}
Let $A$ be a unital \cst-algebra and $e\in \IA$. We define
\begin{equation*}
\I_{c_e}:=\{f\in \IA\colon c_f=c_e\}, \ \ \ \I_{\bar{c_e}}:= \{f\in\IA\colon c_f=c_e^{-1}\}.
\end{equation*}
\end{notation}

\begin{lemma}\label{lem7}
Let $(A,B,\varphi,\theta)$ be as in (\ref{note}) with $A$ oddly decomposable. Let $e,f$ be two non-trivial idempotents in $A$. Then there exist idempotents $e',f'\in \IA$ and $u\in GL(A)$ such that
\begin{equation*}
e',f'\in \I_{c_f}, \ \ \ e'e=ee'=e'\neq e, \ \ \ f'f=ff'=f'\neq f, \ \ \ e'=uf'u^{-1}.
\end{equation*}
\end{lemma}

\begin{proof}
One can adopt the proof of Corollary~2.17 in \cite{AlRBooGio} to the present situation. We include the details for completeness: Find an odd integer $n\geq 3$ and $n$ commuting, orthogonal, idempotents $g_1,\dots,g_n\in \IA$ adding to $f$, such that each $g_i$ is similar to some $g_i'\in\IA$ with the property that $g_i'e=eg_i'=g_i'\neq e$. 

For each $i=1,\dots,n$ set $e_i:=f-g_i$. By Corollary~\ref{cor2} we get that for each $i$ either $c_{g_i}=c_f$ or $c_{g_i}=c_f^{-1}$. If $c_{g_i}=c_f^{-1}$ for $i=1,2,3$ then Lemma~\ref{lem6} ensures that $c_{g_1+g_2+g_3}=c_f^{-1}$. If $c_{g_i}=c_f^{-1}$ for all $i=1,\dots,5$ (if applicable) Lemma~\ref{lem6}, used on $g_1+g_2+g_3$, $g_4$, $g_5$,  ensures that $c_{g_1+\dots+g_5}=c_f^{-1}$. By induction we conclude that $c_{g_1+\dots+g_n}=c_f^{-1}$ if $c_{g_i}=c_f^{-1}$ for all $i$. Knowing that $c_{g_1+\dots+g_n}=c_f\neq c_f^{-1}$ we conclude that $c_{g_m}=c_f$ for some $m\in\{1,\dots,n\}$. Define $f':=g_m$ and $e':=g_m'$. Since $f'$ is similar to $e'$ there exist $u\in GL(A)$ such that $e'=uf'u^{-1}$. Using Proposition~\ref{prop2}$(i)$, $e'\in I_{c_f}$. We conclude that
\begin{equation*}
e',f'\in \I_{c_f}, \ \ \ e'e=ee'=e'\neq e, \ \ \ f'f=ff'=f'\neq f, \ \ \ e'=uf'u^{-1}.\qedhere
\end{equation*}
\end{proof}

\begin{lemma}\label{lem8}
Let $(A,B,\varphi,\theta)$ be as in (\ref{note}) with $A$ oddly decomposable. Then for each $e\in \IA$
\begin{equation*}
\IA=\I_{c_e}\cup \I_{\bar{c_e}}.
\end{equation*}
\end{lemma}

\begin{proof}
The proof of Remark 3.2 \cite{AlRBooGio} does not generalise nicely, so we include a short proof:
Fix any $f\in \IA$. Lemma~\ref{lem7} provides an idempotent $e'\in \IA$ such that
\begin{equation*}
e'\in \I_{c_f}, \ \ \ e'e=ee'=e'\neq e.
\end{equation*}
Applying Corollary~\ref{cor2} to $e'$ and $e-e'$ we get that either $c_{e'}=c_e$ or $c_{e'}=c_e^{-1}$. Hence $c_{f}=c_e$ or $c_{f}=c_e^{-1}$.
\end{proof}

\begin{remark}
Borrowing material from a forthcoming paper \cite{GioSie4} let us mention the following result: Let $(A,B,\varphi,\theta)$ be as in (\ref{note}) with $\varphi$ continuous. Then for each $e\in \IA$
\begin{equation*}
\IA=\I_{c_e}\cup \I_{\bar{c_e}}.
\end{equation*}
\end{remark}

\begin{lemma}\label{lem9}
Let $(A,B,\varphi,\theta)$ be as in (\ref{note}) with $A$   oddly decomposable. Let $e,f$ be two non-trivial orthogonal idempotents in $A$ not adding to one. Suppose that $\theta$ preserves (resp.\ flips) orthogonality on $\{e,f\}$. Then $\theta$ preserves (resp.\ flips) orthogonality on all of $\I_{c_e}$.
\end{lemma}

\begin{proof}
The proof of Lemma~3.4 in \cite{AlRBooGio} can be adapted to the present situation. We include the details for completeness: Let $g,h\in \I_{c_e}$ be any two commuting, orthogonal idempotents. It is enough to show $\theta$ preserves (resp.\ flips) orthogonality on $\{g,h\}$. If $g+h=1$ then $\theta(g)\theta(h)=\theta(g)\theta(1-g)=0$, hence assume $g+h\neq 1$. Using Lemma \ref{lem7} select $x,x',y',z'\in \IA$ and $u\in GL(A)$ such that
\begin{align*}
&x,x'\in \I_{c_{1-e-f}},\ \ \ x{(1-g-h)}={(1-g-h)}x=x\neq {1-g-h},\\
&\ \ \  x'{(1-e-f)}={(1-e-f)}x'=x'\neq {1-e-f},\ \ x=ux'u^{-1}\\
&y'\in \I_{c_{e}},\ \ \  y'x'=x'y'=y'\neq x',\\
&z'\in \I_{c_{e}},\ \ \  z'(x'-y')=(x'-y')z'=z'\neq x'-y'.
\end{align*}
By assumption $\theta$ preserves (resp.\ flips) orthogonality on $\{e,f\}\subseteq \I_{c_e}$.
Using Lemma \ref{lem5} we get that $\theta$ preserves (resp.\ flips) orthogonality on the evidently commuting, orthogonal idempotents 
\begin{equation*}
\{z',y',e,f\}\subseteq \I_{c_e}.
\end{equation*}
Define $y:=uy'u^{-1}$ and $z:=uz'u^{-1}$. Since $y',z'\in \I_{c_e}$ we obtain that $y,z\in \I_{c_e}$, see Proposition \ref{prop2}. By Proposition \ref{prop1}, $\theta(z)\theta(y)=\varphi(u)\theta(z')\theta(y')\varphi(u)^{-1}=0$ (resp. $\theta(1-z)\theta(1-y)=\varphi(u)\theta(1-z')\theta(1-y')\varphi(u)^{-1}=0$). Now Lemma~\ref{lem5} ensures that $\theta$ preserves (resp.\ flips) orthogonality on the clearly commuting, orthogonal idempotents
\begin{equation*}
\{z,y,g,h\}\subseteq \I_{c_e}.\qedhere
\end{equation*}
\end{proof}

\begin{lemma}\label{lem10}
Let $(A,B,\varphi,\theta)$ be as in (\ref{note}) with $A$   oddly decomposable. Let $e,f$ be two non-trivial idempotents that are not $\is$-equivalent. If $\theta$ preserves (resp.\ flips) orthogonality on one of the sets $\I_{c_e}$, $\I_{c_f}$, then $\theta$ flips (resp.\ preserves) orthogonality on the other set.
\end{lemma}

\begin{proof}
One can adopt the proof of Proposition~3.6 in \cite{AlRBooGio} to the present situation. We include the details for completeness: It suffices to show that
\begin{align} 
\textrm{$\theta$ can not preserve orthogonality on $\I_{c_e}\cup\I_{c_f}$.}\label{d3.1}\\
\textrm{$\theta$ can not flip orthogonality on $\I_{c_e}\cup\I_{c_f}$.}\label{d3.2}
\end{align} 

Let us argue why (\ref{d3.1})-(\ref{d3.2}) suffice: Suppose $\theta$ preserves (resp.\ flips) orthogonality on $\I_{c_e}$. Using Lemma~\ref{lem7} select orthogonal idempotents $g,h\in \I_{c_f}$ not adding to one. By Corollary~\ref{cor1} and Lemma~\ref{lem9} we obtain that $\theta$ either preserves or flips orthogonality on all of $\{g,h\}$, and hence on all of $\I_{c_f}$. We conclude $\theta$ flips (resp.\ preserves) orthogonality on $\I_{c_f}$ by (\ref{d3.1})-(\ref{d3.2}).

\eqref{d3.1}: By Lemma~\ref{lem8}, $c_e=c_f^{-1}$. Using Lemma~\ref{lem7} choose $x,y,z\in \IA$ such that
\begin{align*}
&x\in \I_{c_{f}},\ \ \ x{(1-e)}={(1-e)}x=x\neq {1-e},\\
&y\in \I_{c_{f}},\ \ \  ye=ey=y\neq e,\\
&z\in \I_{c_{e}},\ \ \  zx=xz=z\neq x.
\end{align*}
If $1-e-x\in \I_{c_{f}}$ then $\{x,y,1-e-x\}\subseteq \I_{c_{f}}$ are commuting, orthogonal $\is$-equivalent idempotents in $A$ not adding to one. By Lemma~\ref{lem6} we have that $1-e+y\in \I_{c_{f}}$. Hence $e-y\in \I_{c_{f}}$, see Proposition~\ref{prop2}. In particular
\begin{equation*}
c_y=c_{e-y}=c_{e}^{-1}.
\end{equation*}
If $1-e-x\in \I_{c_{e}}$ then $\{e,z,1-e-x\}\subseteq \I_{c_{e}}$ are commuting, orthogonal, $\is$-equivalent idempotents in $A$ not adding to one. It follows that $1-x+z\in \I_{c_{e}}$ and $x-z\in \I_{c_{e}}$. In particular
\begin{equation*}
c_z=c_{x-z}=c_{x}^{-1}.
\end{equation*}
We conclude that $\theta$ can not preserve orthogonality on $\I_{c_{e}}\cup\I_{c_{f}}$.

\eqref{d3.2}: By Lemma~\ref{lem8}, $c_e=c_f^{-1}$. Using Lemma~\ref{lem7} choose $x,y,z\in \IA$ and such that
\begin{align*}
&x\in \I_{c_{f}},\ \ \ x{(1-e)}={(1-e)}x=x\neq {1-e},\\
&y\in \I_{c_{e}},\ \ \  ye=ey=y\neq e,\\
&z\in \I_{c_{f}},\ \ \  zx=xz=z\neq x.
\end{align*}
If $1-e-x\in \I_{c_{f}}$ then $c_{e-y}\neq c_y$ implies that $c_{e-y}=c_y^{-1}=c_f$, see Lemma~\ref{lem8}. Hence $\{e-y,x,1-e-x\}\subseteq \I_{c_{f}}$ are commuting, orthogonal, $\is$-equivalent, idempotents in $A$ not adding to one. It follows that $1-y\in \I_{c_{f}}$ and $y\in \I_{c_{f}}$. Contradiction. In particular
\begin{equation*}
c_y=c_{e-y}=c_{e}.
\end{equation*}
If $1-e-x\in \I_{c_{e}}$ then $c_{x-z}\neq c_z$ implies that $c_{x-z}=c_z^{-1}=c_e$, see Lemma~\ref{lem8}. Hence $\{x-z,e,1-e-x\}\subseteq \I_{c_{e}}$ are commuting, orthogonal, $\is$-equivalent, idempotents in $A$ not adding to one. It follows that $1-z\in \I_{c_{e}}$ and $z\in \I_{c_{e}}$. Contradiction. In particular
\begin{equation*}
c_z=c_{x-z}=c_{x}.
\end{equation*}
We conclude that $\theta$ can not flip orthogonality on $\I_{c_{e}}\cup\I_{c_{f}}$.
\end{proof}

\begin{theorem}\label{thm2}
Let $(A,B,\varphi,\theta)$ be as in (\ref{note}). If $A$ is  oddly decomposable then $\varphi$ induces an orthoisomorphism between the sets of idempotents $\I(A)$ and $\I(B)$, which preserves similarity of idempotents.
\end{theorem}

\begin{proof}
The proof of Theorem~2.21 in \cite{AlRBooGio} can be adapted to the present situation. We include the details for completeness: We may assume $\IA$ is non-empty. Using Lemma~\ref{lem7} select two non-trivial orthogonal, $\is$-equivalent, idempotents $e,f \in \IA$ not adding to one. Define $o:=c_{e+f}$. By Corollary~\ref{cor2} either $c_{e}=o$ or $c_{e}=o^{-1}$. If $c_{e}=o^{-1}$ then $c_e=c_f=c_{e+f}^{-1}$ and $\theta$ flips orthogonality on $\I_{c_e}$. If $c_{e}=o$ then $c_e=c_f=c_{e+f}$ and $\theta$ preserves orthogonality on $\I_{c_e}$. In any case Lemma~\ref{lem10} ensures that $\theta$ preserves orthogonality on $\I_{o}$ and flips orthogonality on $\I_{\bar{o}}$. Define 
\begin{align*}
\tilde\theta(g)=\left\{
\begin{array}{ll}
\theta(g),&\mbox{ if } g \in \I_o\\
1-\theta(g),&\mbox{ if } g \in \I_{\bar o}\\
1,&\mbox{ if } g =1\\
0,&\mbox{ if } g =0
\end{array}\right.
\end{align*}
Fix two commuting, orthogonal idempotents $g,h\in \I(A)$. If $g$ or $h$ is equal to zero then obviously $\tilde\theta(g)\tilde\theta(h)=0$. If $h,g$ are both nonzero and add to one then $c_h=c_g$, see Proposition~\ref{prop2}. Hence $\tilde\theta$ restricts to either $\theta$ or $1-\theta$ on $\{g,h\}$, implying $\tilde\theta(g)\tilde\theta(h)=0$. We may assume $g,h\in\IA$ with $g+h\neq 1$. Then
\begin{align*}
g,h\in \I_o && \theta(g)\theta(h)=0&& \\
g,h\in \I_{\bar{o}}& \Rightarrow &\theta(1-g)\theta(1-h)=0 &\Rightarrow & \tilde\theta(g)\tilde\theta(h)=0 \\
g\in \I_o, h\in \I_{\bar{o}}&& \theta(g)\theta(1-h)=0&&
\end{align*}
The fact that $\theta(g)\theta(1-h)=0$ above follows indirectly: Since $c_g\neq c_h$ either $\theta(1-g)\theta(h)=0$ or $\theta(g)\theta(1-h)=0$ (by Theorem~\ref{thm1}). The first equality implies that $c_g=c_h^{-1}=c_{g+h}$, hence $g,1-g-h\in \I_o$. But then $\theta(g)\theta(1-g-h)=0$, meaning that $c_g=c_{1-g-h}=c_{1-h}$. Contradiction. We conclude that 
\begin{equation*}
\tilde\theta\colon \I(A)\to \I(B)
\end{equation*}
preserves orthogonally on $\I(A)$. Surjectivity of $\tilde\theta$ follows from Lemma~\ref{lem3}. 
Injectivity of $\tilde\theta$ follows from the fact that if $g\in \I_o$, $h\in \I_{\bar{o}}$ then $\tilde\theta(g)\neq \tilde\theta(h)$, since $\tilde\theta(g)=\tilde\theta(h)$ implies $g=1-h$. Finally, for any $g\in \IA$ and $u\in GL(A)$, we have $g\is ugu^{-1}$, by Proposition~\ref{prop2}, so
\begin{align*}
g\in\I_o&&\tilde\theta(ugu^{-1})&=\theta(ugu^{-1})=\varphi(u)\theta(g)\varphi(u)^{-1}=\varphi(u)\tilde\theta(g)\varphi(u)^{-1},\\
g\in\I_{\bar o}&&\tilde\theta(ugu^{-1})&=1-\theta(ugu^{-1})=\dots=\varphi(u)\tilde\theta(g)\varphi(u)^{-1}.
\end{align*}
We conclude that $\theta$ preserves similarity of idempotents.
\end{proof}

\section{The case of simple AH-algebras}
\subsection{From an orthoisomorphism to a $K_0$-order isomorphism}
\label{sec4.1}
In this subsection, we prove that an (abstract) isomorphism $GL(A)\cong GL(B)$ between the general linear groups of certain stably finite \cst-algebras $A$ and $B$ of real rank zero (including the simple AH-algebras of slow dimension growth) induces an isomorphism between their ordered $K_0$-groups. In particular, we have that if $A$ and $B$ are either two simple unital AF-algebras, or two irrational rotation algebras, then $A$ is $^*$-isomorphic to $B$ if and only if their general linear groups are isomorphic (as abstract groups). Our approach uses ideas of \cite{AlRBooGio}, but with proofs that are somehow different, and some clarifications are given.

\begin{notation}\label{note2}
(i) Let $\F$ denote the class of simple, unital, separable \cst-algebras of real rank zero with cancellation (or equivalently with stable rank one, see Corollary 6.5.7 in \cite{Bla}). Let $\F_1$ denote the class of \cst-algebras $A$ in $\F$ for which $K_0(A)$ is noncyclic and weakly unperforated.

(ii) Let $A$ be a unital \cst-algebra. Denote by $\cP(A)$ the set of projections in $A$, and by $\PA$ the set $\cP(A)\backslash\{0,1\}$ of \emph{non-trivial} projections in $A$.
\end{notation}

The following lemma is well known, see Proposition~4.6.2 and Proposition~4.6.5 in \cite{Bla}.
\begin{lemma}\label{lem11}
Let $A$ be a unital \cst-algebra. Every idempotent in $A$ is similar to a projection in $A$. Every pair of projections is $A$ are similar if, and only if, they are unitary equivalent.
\end{lemma}

\begin{proposition}\label{prop3}
Each \cst-algebra in $\F_1$ is  oddly decomposable.
\end{proposition}
\begin{proof}
The result follows immediately from Proposition~4.2 in \cite{AlRBooGio} in combination with Remark \ref{rem3.1}.
\end{proof}

Following \cite{Bla} and \cite{RorLarLau} an {\em ordered (abelian) group} $G$ is an abelian group with a distinguished {\em positive cone}, \ie{}a subset $G_+\subseteq G$ fulfilling that
\begin{equation*}
G_+ + G_+ \subseteq G_+, \ \ \  G_+ \cap (-G_+)=\{0\}, \ \ \ G_+ - G_+ = G.
\end{equation*}
The set $G_+$ induces a translation-invariant partial ordering on $G$ by $x\leq y$ if $y-x\in G_+$. 

Essentially\footnote{Effros presumes the group is unperforated (resp.\ is a dimension group) in his definition of an ordered group (resp.\ a scaled dimension group). We have removed these two constraints and changed the terminology accordingly.} as in \cite{Eff}, a \emph{scaled ordered group} $G$ is a ordered group with a distinguished \emph{scale}, \ie{}a subset $\Gamma = \Gamma (G)$ of $G_+$, which is generating, hereditary and directed, \ie{}

\begin{itemize}
\item[(i)] For each $a \in G_+$, there exist $a_1,\dots , a_r\in \Gamma$ with $a=\sum_{i=1}^r a_i$.
\item[(ii)] If $0\leq a\leq b \in \Gamma$, then $a\in \Gamma$.
\item[(iii)] Given $a,b\in \Gamma$, there exist $c\in\Gamma$ with $a,b\leq c$.
\end{itemize}

A scale $\Gamma$ has a partially defined addition; in fact $a\geq b$ in $\Gamma$ if, and only if $a=b+c$ for some $c\in \Gamma$. Following \cite{Eff}, a group homomorphism of scaled ordered groups $\alpha\colon G \to G'$ is a \emph{contraction} if $\alpha(\Gamma(G)) \subseteq \Gamma(G')$. If $\Gamma$ and $\Gamma'$ are scales of two scaled ordered groups, then a map $\alpha\colon \Gamma \to \Gamma'$ is a \emph{scale homomorphism} (resp.\ a \emph{scale isomorphism}) if $a = b + c$ in $\Gamma$ implies that (resp.\ is equivalent to) $\alpha(a) = \alpha(b) + \alpha(c)$ in $\Gamma'$. 

\begin{proposition}[Effros]\label{prop4}
Let $G$ and $G'$ be two scaled ordered groups with Riesz interpolation. Any scale homomorphism $\alpha\colon \Gamma (G) \to \Gamma (G')$ extends to a unique contraction $\tilde{\alpha}\colon G \to G'$. If $\alpha$ is a scale isomorphism, then $\tilde{\alpha}$ is an isomorphism of the scaled ordered groups $G$ and $G'$.
\end{proposition}
\begin{proof}
Notice that the proof of Lemma~7.3 and Corollary~7.4 in \cite{Eff} does not use perforation nor countability of the groups involved.
\end{proof}	

\begin{lemma}\label{lem4.4}
If $A\in \F$, then $K_0(A)$ is a simple scaled ordered group with Riesz interpolation and scale $\Sigma(A):=\{[p]\colon p\in \cP(A)\}$.
\end{lemma}

\begin{proof}
Since $A$ is stably finite and simple, the group $K_0(A)$ is a simple scaled ordered group with Riesz interpolation by Proposition~3.3.7 and Theorem~3.3.18 in \cite{Lin}, and $\Sigma(A)$ is hereditary and directed, see p.\ 38 in \cite{Bla}.

For sake of completeness we show $\Sigma(A)$ is generating: Fix any $x$ in $K_0(A)_+$. Recall that $x=[p]$ for some projection $p$ in $M_n(A)$ (with $n\in \NN$). Let $1_n$ denote the unit of $M_n(A)$. Cleary, $1_n=\sum_{i=1}^n e_{ii}$, where $e_{ii}\in M_n(A)$ is the matrix with 1 at entry $(i,i)$ and zero otherwise, and $p\leq 1_n$. Since $M_n(A)$ has real rank zero it follows from Corollary~3.3.17 in \cite{Lin} that there exist projections $p_i\in M_n(A)$ such that $[p_i]\leq [e_{ii}]$ and $\sum_{i=1}^n p_{i}=p$. Hence $x=\sum_{i=1}^n [p_{i}]$ and $[p_{i}]\in \Sigma(A)$, using the characterisation $\Sigma(A)=\{x\in K_0(A)_+\colon x\leq [1]\}$, see \cite{Bla}.
\end{proof}

\begin{theorem}
\label{K0isom}
Let A and B be two \cst-algebras in $\F_1$. If $GL(A)$ and $GL(B)$ are isomorphic (as abstract groups), then $K_0(A)$ and $K_0(B)$ are isomorphic as scaled ordered groups.
\end{theorem}

\begin{proof}
Let $\tilde\theta \colon \I(A) \to \I(B)$ be the orthoisomorphism preserving similarity of idempotents given by Theorem~\ref{thm2} and Proposition~\ref{prop3} with $\tilde\theta(1)=1$. Let 
\begin{equation*}
\tilde\theta_*\colon \Sigma(A)\to \Sigma(B)
\end{equation*} be given by $\tilde\theta_*([p])=[p']$, where $p'=u\tilde\theta(p)u^{-1}$ for some $u\in GL(B)$ such that $p'$ is a projection in $B$, see Lemma~\ref{lem11}.

We show $\tilde\theta_*$ is well defined: Fix any two projections $p,q$ in $\cP(A)$. Assume $[p]=[q]$. By Proposition~3.1.7(iv) in \cite{RorLarLau} also $[1-p]=[1-q]$. Since $A$ has cancellation $p$ and $q$ are unitary equivalent, see Definition~7.3.1 and Proposition~2.2.2 in \cite{RorLarLau}. Since $\tilde\theta$ preserves similarity $\tilde\theta(p)$ and $\tilde\theta(q)$ are similar. We conclude that $p'$ and $q'$ are similar and (by Lemma~\ref{lem11}) unitary equivalent. We obtain that $[p']=[q']$.

We show $\tilde\theta_*$ is a scale homomorphism: Fix $x,y,z\in \Sigma(A)$ satisfying the equality $x+y=z$. Find $p,q\in \cP(A)$ such that $x=[p]$ and $y=[q]$. Select $q_1\in \cP(A)$ such that 
\begin{equation*}
[q_1]=[p]\ \ \textrm{ and } \ \ q_1\leq 1-q
\end{equation*}
as follows: If $z=[1]$ set $q_1:=1-q$ implying $[q_1]=z-y=[p]$. If $z<[1]$ use $[p]=z-[q]<[1-q]$ and Corollary~6.9.2 in \cite{Bla} to deduce that $p$ is Murray-von Neumann equivalent to a subprojection, say $q_1$, of $1-q$. (Here we have used $A$ is simple, unital, stably finite, of real rank zero, with cancellation and weakly unperforated $K_0(A)$.) We obtain that 
\begin{equation*}
\tilde\theta(q_1)+\tilde\theta(q)=\tilde\theta(q_1+q)\ \ \textrm{ and } \ \ q_1q=0.
\end{equation*}
Find $v\in GL(B)$ such that $v\tilde\theta(q_1)v^{-1}$ and $v\tilde\theta(q)v^{-1}$ are orthogonal projections in $B$ (easy exercise, see Remark \ref{rem3.1}). Hence
\begin{align*}
\tilde\theta_*(z)&= \tilde\theta_*([p]+[q])=\tilde\theta_*([q_1]+[q])=\tilde\theta_*([q_1+q])\\
&=[v\tilde\theta(q_1+q)v^{-1}]=[v\tilde\theta(q_1)v^{-1}]+[v\tilde\theta(q)v^{-1}]\\
&=\tilde\theta_*([q_1])+\tilde\theta_*([q])=\tilde\theta_*(x)+\tilde\theta_*(y)
\end{align*}

We show $\tilde\theta_*$ is a scale isomorphism: As $\tilde\theta : \I(A) \to \I(B)$ is an orthoisomorphism, its inverse induces a scale homomorphism $(\tilde\theta^{-1})_*\colon \Sigma(B)\to \Sigma(A)$.
For $p\in \cP(A)$ we have that 
\begin{align*}
(\tilde\theta^{-1})_*( \tilde\theta_*([p]))&=(\tilde\theta^{-1})_*[u\tilde\theta(p)u^{-1}]
=[v\tilde\theta^{-1}(u\tilde\theta(p)u^{-1})v^{-1}]\\
&=[vwpw^{-1}v^{-1}]
=[p],
\end{align*}
for appropriate $u\in GL(B)$, and $v,w\in GL(A)$, using that $\tilde\theta^{-1}$ maps $u\tilde\theta(p)u^{-1}$ to an idempotent similar to $p$. By symmetry both $(\tilde\theta^{-1})_*\circ \tilde\theta_*$ and $\tilde\theta_*\circ(\tilde\theta^{-1})_*$ are identity maps. Hence $(\tilde\theta^{-1})_*=(\tilde\theta)_*^{-1}$. 

Using Proposition~\ref{prop4} and Lemma~\ref{lem4.4} we obtain that $K_0(A)$ and $K_0(B)$ are isomorphic as scaled ordered groups.
\end{proof}

\begin{lemma}
\label{lemm5}
Every infinite-dimensional, simple, unital AH-algebras of slow dimension growth and of real rank zero belongs to the class $\F_1$.
\end{lemma}
\begin{proof}
Cancellation: We refer to Theorem~1 in \cite{BlaDadRor} and Proposition~6.5.1 in \cite{Bla}. Weakly unperforated: See p.\ 2 in \cite{Vil}. Noncyclic: See Remark 2.7 in \cite{Dad}.
\end{proof}	

\begin{corollary}\label{cor3}
If $A$ and $B$ are simple, unital AH-algebras of slow dimension growth and of real rank zero, with isomorphic general linear groups (as abstract groups), 
then 
\begin{equation*}
(K_0(A), K_0(A)_+ , [1_A])\ \  \mbox{ and }\ \ (K_0(B),K_0(B)_+,[1_B])
\end{equation*}
are order isomorphic by a map preserving the distinguished order units.
\end{corollary}
\begin{proof}
If $A$ is infinite-dimensional then so is $B$. (By Lemma~\ref{lemm5} and Proposition~\ref{prop3} $A$ is  oddly decomposable. We can therefore find arbitrary many  orthogonal idempotents $(e_i)$ in $\IA$. The isomorphism of $GL(A)$ and $GL(B)$ induces an orthoisomorphism $\tilde\theta \colon \I(A) \to \I(B)$, by Theorem~\ref{thm2}. The orthogonal idempotents $\tilde\theta(e_i)$ in $\IB$ ensure $B$ is infinite-dimensional.) The desired result follows now from Theorem~\ref{K0isom}. If both $A$ and $B$ are finite dimensional we refer to \cite{GioSie3}.
\end{proof}
	
Using H. Lin's characterization of TAF-algebras (see \cite{Lin} or Theorem~3.3.5 in \cite{RorSto}) we can also state Corollary~\ref{cor3} as follows.

\begin{corollary}\label{cor4.8}
Let $A$ and $B$ be two simple, unital, nuclear, separable TAF-algebras of real rank zero, in the UCT-class $\N$ with isomorphic general linear groups (as abstract groups), then 
\begin{equation*}
(K_0(A), K_0(A)_+ , [1_A])\ \  \mbox{ and }\ \ (K_0(B),K_0(B)_+,[1_B])
\end{equation*}
are order isomorphic by a map preserving the distinguished order units.
\end{corollary}

\begin{corollary}\label{cor4.9}
If $A$ and $B$ are either two simple unital AF-algebras, or two irrational rotation algebras, then $A$ is $^*$-isomorphic to $B$ if and only if their general linear groups are isomorphic (as abstract groups).
\end{corollary}
\begin{proof}
Both the class of unital simple AF-algebras and the class of irrational rotation algebras are classified by $(K_0, {K_0}_+, [1])$, see Theorem~7.3.4 in \cite{RorLarLau} and Corollary~VI.5.3 in \cite{Dav}. 

Any unital simple AF-algebras is a nuclear TAF-algebras of real rank zero, in the UCT-class $\N$, and any irrational rotation algebra is an AH-algebras of slow dimension growth and of real rank zero, see \cite{EllEva}.
\end{proof}

\subsection{From a general linear group isomorphism to a \cst-isomorphism}
\label{sec4.2}

For simple AH-algebras of real rank zero, let us recall the classification theorem, provided independently by Gong in \cite{Gon} and Dadarlat in \cite{Dad2}, whose proof uses Elliott-Gong's classification in \cite{EllGon} (see for example Theorem~3.3.1 in \cite{RorSto}).

\begin{theorem}[Dadarlat, Gong, Elliott]\label{thm15}
Let $A$ and $B$ be simple, unital, AH-algebras of slow dimension growth and of real rank zero. It follows that $A$ is $^*$-isomorphic to $B$ if and only if
\begin{equation*}
(K_0(A), K_0(A)_+ , [1_A])\cong (K_0(B), K_0(B)_+ , [1_B]), \ \ \ K_1(A)\cong K_1(B).
\end{equation*}
\end{theorem}

\begin{notation}\label{note6}
Let $A$ be a unital \cst-algebra. We equip the general linear group $GL(A)$ with the topology induced by the norm on $A$. Denote by $GL_0 (A)$ the connected component $\{u\colon u\hsim 1\}$ of the identity element in $GL(A)$.
\end{notation}

\begin{theorem}\label{thm6}
Let $A$ and $B$ are simple, unital AH-algebras of slow dimension growth and of real rank zero. Then $A$ and $B$ are isomorphic if and only if their general linear groups are topologically isomorphic.
\end{theorem}
\begin{proof}
If $A$ and $B$ are isomorphic then their general linear groups are topologically isomorphic. Conversely 
let $\varphi \colon GL(A) \to GL(B)$ be a topological isomorphism from $GL(A)$ onto $GL(B)$. By continuity of $\varphi$, $\varphi(GL(A)_0)=GL(B)_0$. It follows that $u+GL(A)_0\mapsto \varphi(u)+GL(B)_0$ is an isomorphism of $GL(A)/GL(A)_0$ and $GL(B)/GL(B)_0$ (with inverse $v+GL(B)_0\mapsto \varphi^{-1}(v)+GL(A)_0$). Recall that for a unital \cst-algebra $C$ of stable rank one, $K_1(C)$ is isomorphic to the group $GL(C)/GL(C)_0$, see Theorem~2.10 in \cite{Rie}. Consequently, we conclude that $K_1(A)$ is isomorphic to $K_1(B)$.
\end{proof}

\section{The case of Kirchberg algebras}
\label{sec5}
\subsection{From an orthoisomorphism to a $K_0$-isomorphism}
In this subsection, inspired by \cite{AlRBooGio}, we show that an isomorphism between the general linear groups of simple, unital, purely infinite \cst-algebras induces an isomorphism between their $K_0$-groups.

\begin{theorem}\label{thm5.1}
Every simple, unital, purely infinite \cst-algebra is  oddly decomposable.
\end{theorem}

\begin{proof}
The result follows immediately from Theorem~5.2 in \cite{AlRBooGio} in combination with Remark \ref{rem3.1}.
\end{proof}

Recall that if $A$ is a purely infinite simple \cst-algebra, then every nonzero projection in $A$ is infinite, and $K_0(A)=\{[p]\colon p\in \cP(A), p\neq 0\}$ (see p.\ 73--85 in \cite{RorLarLau}). If $A$, in addition, is unital then $1$ is an infinite projection and therefore Murray-von Neumann equivalent to a subprojection $q < 1$. Hence $[q]=[1]$, and $K_0(A)=\{[p]\colon p\in \PA\}$.

\begin{theorem}\label{thm3}
If $A$ and $B$ are two unital, simple, purely infinite \cst-algebras, whose general linear groups are isomorphic (as abstract groups), then there is an isomorphism from $K_0(A)$ to $K_0(B)$, sending $[1_A]$ to $[1_B]$.
\end{theorem}

\begin{proof}
Let $\tilde\theta \colon \I(A) \to \I(B)$ be the orthoisomorphism preserving similarity of idempotents given by Theorem~\ref{thm2} and Theorem~\ref{thm5.1} with $\tilde\theta(1)=1$. Let 
\begin{equation*}
\tilde\theta_*\colon K_0(A)\to K_0(B)
\end{equation*}
be given by $\tilde\theta_*([p])=[p']$, where $p'=u\tilde\theta(p)u^{-1}$ for some $u\in GL(B)$ such that $p'$ is a projection in $B$, see Lemma~\ref{lem11}.

We show $\tilde\theta_*$ is well defined and injective: Fix any two projections $p,q$ in $\PA$. Assume $[p]=[q]$. Since $p,q$ are infinite the assumption is equivalent to $p,q$ being unitary equivalent, see Corollary~6.11.9 in \cite{Bla}. By Lemma~\ref{lem11} the assumption is equivalent to $p,q$ being similar. Since $\tilde\theta$ and $\tilde\theta^{-1}$ preserves similarity the assumption is equivalent to $\tilde\theta(p),\tilde\theta(q)$ being similar. By definition of $p',q'$ the assumption is equivalent to $p',q'$ being similar, and hence unitary equivalent (see Lemma~\ref{lem11}). Using Corollary~6.11.9 in \cite{Bla} once more we obtain that the assumption $[p]=[q]$ is equivalent to $[p']=[q']$.

We show $\tilde\theta_*$ is unital and surjective: Since $\tilde\theta(1)=1$ is a projection we get that $\tilde\theta_*([1])=[\tilde\theta(1)]=[1]$. Fix a projection $p\in \PB$. Find an idempotent $e\in \IA$ such that $\tilde\theta(e)=p$. By Lemma~\ref{lem11} there exist a $u\in GL(A)$ such that $ueu^{-1}$ is a projection. Now 
\begin{equation*}
\tilde\theta_*([ueu^{-1}])=[v\tilde\theta(ueu^{-1})v^{-1}],
\end{equation*}
for an appropriate $v\in GL(B)$. Since $e$ is similar to $ueu^{-1}$ then $\tilde\theta(e)$ is similar to $v\tilde\theta(ueu^{-1})v^{-1}$. By Lemma~\ref{lem11} similar projections are unitary equivalent. Hence Corollary~6.11.9 in \cite{Bla} ensures $\tilde\theta_*([ueu^{-1}])=[p]$.

We show $\tilde\theta_*$ is a homomorphism: Fix any $p,q\in \PA$. Since $1-q,q$ are (full and properly) infinite we can find projections $r\leq 1-q$ and $s\leq q$ in $A$ such that $p$ (resp.\ $q$) is Murray-von Neumann equivalent to $r$ (resp.\ $s$), see p.\ 75 in \cite{RorLarLau}. In particular $[p]=[r]$ and $[q]=[s]$. Since $r$ and $s$ are orthogonal $[r+s]=[r]+[s]$. Find $v\in GL(B)$ such that $v\tilde\theta(r)v^{-1}$ and $v\tilde\theta(s)v^{-1}$ are orthogonal projections in $B$ (see Remark \ref{rem3.1}). Hence
\begin{align*}
\tilde\theta_*([p]+[q])&=\tilde\theta_*([r]+[s])=\tilde\theta_*([r+s])=[v\tilde\theta(r+s)v^{-1}]\\
&=[v\tilde\theta(r)v^{-1}+v\tilde\theta(s)v^{-1}]=[v\tilde\theta(r)v^{-1}]+[v\tilde\theta(s)v^{-1}]\\
&=\tilde\theta_*([r])+\tilde\theta_*([s])=\tilde\theta_*([p])+\tilde\theta_*([q])
\end{align*}
This shows that $\tilde\theta_*$ is the desired isomorphism.
\end{proof}

In \cite{Cun}, J. Cuntz proved that for $2\leq n<\infty$, $K_0(\cO_n) \cong  \ZZ/(n - 1)\ZZ$ and $K_0(\cO_\infty) \cong \ZZ$. Hence, we have:

\begin{corollary}
\label{cor5.3}
Two Cuntz algebras are isomorphic if and only if their general linear groups are isomorphic (as abstract groups).
\end{corollary}

\subsection{From a general linear group isomorphism to a \cst-isomorphism}
\label{sec5.2}
Recall that a \emph{Kirchberg algebra} is a purely infinite, simple, nuclear, separable \cst-algebra, see Definition~4.3.1 in \cite{RorSto}, and that the following result of Kirchberg and Phillips essentially classifies such algebras.

\begin{theorem}[Kirchberg, Phillips]\label{thm16}
Let $A$ and $B$ be unital Kirchberg algebras in the UCT-class $\N$. Then $A$ and $B$ are $^*$-isomorphic if, and only if, there exist isomorphisms $\alpha_0 \colon K_0(A) \to K_0(B)$ and $\alpha_1\colon K_1(A) \to K_1(B)$ with $\alpha_0([1_A]) = [1_B]$.
\end{theorem}

\begin{notation}\label{note7}
Let $A$ be a unital \cst-algebra. As usual, the topology on the unitary group $\U(A)$ is inherited from $GL(A)$. Denote by $\U_0(A)$ the connected component $\{u\colon u\hsim 1\}$ of the identity element in $\U(A)$.
\end{notation}

\begin{theorem}\label{thm4}
If $A$ and $B$ are two unital, simple, purely infinite \cst-algebras, whose general linear groups are isomorphic (as abstract groups), then the groups $K_1(A)$ and $K_1(B)$ are isomorphic.
\end{theorem}

\begin{proof}
Let $\varphi \colon GL(A) \to GL(B)$ denote the isomorphism from $GL(A)$ into $GL(B)$. Since $\varphi$ preserves symmetries (\ie{}if $s^2=1$ in $GL(A)$ then $\varphi(s)^2=1$ in $GL(B)$) and symmetries generate the connected component of the identity (by Theorem~3.7 in \cite{Lee}) we have that $\varphi(GL_0 (A))=GL_0 (B)$. Let 
\begin{equation*}
\tilde\varphi_*\colon \U(A)/\U_0(A)\to \U(B)/\U_0(B)
\end{equation*} 
be given by $\tilde\varphi_*([u])=[u']$, where $u'=\omega(\varphi(u))$ and $\omega$ is the map from Proposition~2.1.8 in \cite{RorLarLau} turning invertible elements into unitaries.

We show $\tilde\varphi_*$ is well defined and injective: Fix any two unitaries  $u,v$ in $\U(A)$. Assume $[u]=[v]$. Recall that $u\hsim v$ in $\U(A)$ if, and only if $u\hsim v$ in $GL(A)$, see Proposition~2.1.8 in \cite{RorLarLau}. In particular the assumption is equivalent to $\varphi(u)\hsim \varphi(v)$ in $GL(B)$ (recalling $\varphi(GL_0 (A))=GL_0 (B)$). By Proposition~2.1.8 in \cite{RorLarLau} both $\varphi(u)\hsim \omega(\varphi(u))$ and $\varphi(v)\hsim \omega(\varphi(v))$ in $GL(B)$. Hence the assumption is equivalent to  $u'\hsim v'$ in $GL(B)$, and hence also to $[u']=[v']$.

We show $\tilde\varphi_*$ is surjective: Fix an unitary $u\in \U(B)$. Find an invertible element $v\in GL(A)$ such that $\varphi(v)=u$. Similarly to a previous argument we have that $\omega(v)\hsim v= \varphi^{-1}(u)$ in $GL(A)$ and $\varphi(\omega(v))\hsim \varphi(v)= u$ in $GL(B)$. Using that $\varphi(\omega(v))\hsim \omega(\varphi(\omega(v)))$ we obtain that $[\omega(\varphi(\omega(v)))]=[u]$. Hence $\tilde\varphi_*([\omega(v)])=[u]$.

We show $\tilde\varphi_*$ is a homomorphism: Fix any $u,v\in \U(A)$. Using the equivalences $\varphi(u)\hsim \omega(\varphi(u))$ and $\varphi(v)\hsim \omega(\varphi(v))$ in $GL(B)$ we obtain that 
\begin{equation*}
\varphi(u)\varphi(v)\hsim \varphi(u)\omega(\varphi(v)) \hsim \omega(\varphi(u))\omega(\varphi(v))\ \ \ \ \mbox{ in }GL(B).
\end{equation*}
We also have that $\omega(\varphi(uv))\hsim\varphi(uv)$ in $GL(B)$. Combining these relations we have that $[\omega(\varphi(uv))]=[\omega(\varphi(u))\omega(\varphi(v))]$. We conclude that
\begin{equation*}
\tilde\varphi_*([u])\tilde\varphi_*([v])=[\omega(\varphi(u))\omega(\varphi(v))]=[\omega(\varphi(uv))]=\tilde\varphi_*([uv]).
\end{equation*}

This shows that $\tilde\varphi_*$ is the desired isomorphism. Recall that for a unital purely infinite simple \cst-algebra $C$, $K_1(C)$ is isomorphic to $\U(C)/\U_0(C)$ by Theorem~1.9 in \cite{Cun}. Consequently, we conclude that $K_1(A)$ is isomorphic to $K_1(B)$.
\end{proof}

Thanks to Theorems \ref{thm3} and \ref{thm4}, we have the following conclusion:

\begin{corollary}\label{cor5.6}
Let $A$ and $B$ be two unital Kirchberg algebras in the UCT-class $\N$. Then $A$ and $B$ are isomorphic if and only if their general linear groups are isomorphic (as abstract groups).
\end{corollary}

\bibliographystyle{amsplain}
\providecommand{\bysame}{\leavevmode\hbox to3em{\hrulefill}\thinspace}
\providecommand{\MR}{\relax\ifhmode\unskip\space\fi MR }
\providecommand{\MRhref}[2]{%
  \href{http://www.ams.org/mathscinet-getitem?mr=#1}{#2}
}
\providecommand{\href}[2]{#2}

\end{document}